\documentclass[12pt]{amsart}
\usepackage{latexsym,amssymb,amsmath,amsthm,amscd,graphicx}
\usepackage{esint} 

\usepackage{color}
\numberwithin{equation}{section}
\newtheorem{theorem}{Theorem}[section]
\newtheorem{lemma}[theorem]{Lemma}
\newtheorem{corollary}[theorem]{Corollary}
\newtheorem{proposition}[theorem]{Proposition}

\theoremstyle{definition}

\newtheorem{remark}[theorem]{Remark}
\newtheorem{definition}[theorem]{Definition}

\theoremstyle{remark}

\usepackage{url,hyperref}

\title[Characterization for Campanato norm via sBfs]{Characterization for Campanato norm via quasi-Banach function spaces not assuming
\\
the Fatou property}
\author{Naoya Hatano}
\date{}

\address{Graduate School of Information Science and Technology, the University of Osaka, 1-5, Yamadaoka, Suita-shi, Osaka 565-0871, Japan}

\email{n.hatano.chuo@gmail.com}

\begin{document}

\maketitle

\begin{abstract}
It is well known that the BMO and Campanato norms can be characterized using the $L^p$-average.
These characterizations were later generalized to averages taken over various types of function spaces.
In particular, generalizations using Banach function spaces were provided by Ho, Izuki, Noi, and Sawano.
In this paper, as a further generalization, we provide similar characterizations using quasi-Banach function spaces that do not assume the Fatou property.
Note that the duality argument is not available in this setting.
\end{abstract}

{\bf 2010 Classification} 42B35, 46E30


\section{Introduction}

Throughout this paper, let $L^0({\mathbb R}^n)$ denote the space of all measurable functions, and $L_{\rm loc}^1({\mathbb R}^n)$ the space of functions integrable on every compact subset of ${\mathbb R}^n$.
Campanato spaces are generalizations of the John--Nirenberg space ${\rm BMO}({\mathbb R}^n)$.
The space ${\rm BMO}({\mathbb R}^n)$ is defined as the space of all $f\in L_{\rm loc}^1({\mathbb R}^n)$ with bounded mean oscillation, that is,
\[
\|f\|_{\rm BMO}
\equiv
\sup_{Q\in{\mathcal Q}}
\frac1{|Q|}
\int_Q|f(x)-f_Q|\,{\rm d}x
<\infty,
\]
where ${\mathcal Q}({\mathbb R}^n)$ denotes the collection of all cubes in ${\mathbb R}^n$ whose sides are parallel to the coordinate axes, and $f_Q$ is the average of $f$ over $Q$.
It is well known that for $1\le p<\infty$, the BMO norm can be characterized using $L^p$-averages.
More precisely, there exists a constant $C\ge1$ such that for all $f\in L_{\rm loc}^1({\mathbb R}^n)$,
\[
\frac1C
\|f\|_{\rm BMO}
\le
\sup_{Q\in{\mathcal Q}}
\frac{\|(f-f_Q)\chi_Q\|_{L^p}}{\|\chi_Q\|_{L^p}}
\le
C
\|f\|_{\rm BMO},
\]
where $\chi_Q$ is the characteristic function of $Q$.
This characterization has been extended to other types of function spaces by many authors, including weighted Lebesgue spaces \cite{Ho11}, variable exponent Lebesgue spaces \cite{Ho14,IzSa12,IST14}, and Morrey spaces \cite{Ho17,NaSa17}.
Moreover, these results were extended to Banach function spaces in the works of Ho, Isuki, Noi and Sawano (see \cite{Ho12, Ho09, Ho16, Izuki15, INS19}).

Meanwhile, a characterization of the Campanato norm via Lebesgue spaces was given by Nakai \cite{Nakai08}.
The Campanato norm $\|\cdot\|_{{\mathcal L}_{p,\phi}}$ for $1\le p<\infty$ is defined by
\[
\|f\|_{{\mathcal L}_{p,\phi}}
\equiv
\sup_{Q\in{\mathcal Q}}
\phi(Q)
\frac{\|(f-f_Q)\chi_Q\|_{L^p}}{\|\chi_Q\|_{L^p}},
\]
for all locally integrable functions $f$, where $\phi:{\mathbb R}^n\times(0,\infty)\to(0,\infty)$ is a positive function.
For a cube $Q=Q(x,r)=x+[-r,r)^n$, we write $\phi(Q)\equiv\phi(x,r)$.
When $\phi(Q)\equiv1$, the Campanato space ${\mathcal L}_{1,\phi}({\mathbb R}^n)$ coincides with the BMO space.
When $\phi(Q)=\ell(Q)^{-\theta}$ with $0<\theta<1$, where $\ell(Q)$ denotes the side length of $Q$, the Campanato space ${\mathcal L}_{1,\phi}({\mathbb R}^n)$ coincides with the Lipschitz space ${\rm Lip}_\theta({\mathbb R}^n)$ of order $\theta$.
Campanato spaces defined by a function $\phi$ depending on the position of cubes were introduced by Nakai \cite{Nakai06} (see also \cite{Nakai93}).
Lerner \cite{Lerner05} proved that, if $p(\cdot),\,1/p(\cdot)\in L^\infty({\mathbb R}^n)$ and $p(\cdot)\in\mathcal L_{1,\phi}({\mathbb R}^n)$, where
\[
\phi(x,r)
=
\log\left(
e+\max\left(2r,\frac1{2r},|x|\right)
\right),
\quad
(x,r)\in\mathbb{R}^n\times(0,\infty),
\]
then there exists $\alpha>0$ such that the Hardy--Littlewood maximal operator $M$ is bounded on the variable Lebesgue space
$L^{\alpha+p(\cdot)}({\mathbb R}^n)$.
Suppose that $\phi$ is almost decreasing, that is, there exists a constant $C>0$ such that
\begin{equation}\label{eq:phi-dec}
\phi(x,r)
\le C
\phi(x, s)
\quad \text{whenever} \quad
r\ge s,
\end{equation}
and that there exists $C\ge1$ such that for all $x,y\in{\mathbb R}^n$ and $r>0$,
\begin{equation}\label{eq:near}
\frac1C
\le
\frac{\phi(x,r)}{\phi(y,r)}
\le C
\quad \text{whenever} \quad
|x-y|\le r,
\end{equation}
then there exists $C\ge1$ such that
\[
\frac1C
\|f\|_{{\mathcal L}_{1,\phi}}
\le
\|f\|_{{\mathcal L}_{p,\phi}}
\le
C
\|f\|_{{\mathcal L}_{1,\phi}}.
\]

Furthermore, the above characterizations have been generalized using (ball) Banach function spaces by many authors.
For the definition and details of Banach function spaces, we refer the reader to the monograph \cite{BeSh88}.

\begin{definition}\label{def:Bfs}
Let $X({\mathbb R}^n)=(X({\mathbb R}^n),\|\cdot\|_X)\subset L^0({\mathbb R}^n)$ be a quasi-Banach space.
The space $X({\mathbb R}^n)$ is said to be a quasi-Banach function space if it satisfies the following conditions for all $f,g,f_k\in L^0({\mathbb R}^n)$ $(k\in{\mathbb N})${\rm :}
\begin{itemize}
\item[{\rm (i)}] If $f\in X({\mathbb R}^n)$ and $g\in L^0({\mathbb R}^n)$ with $|g|\le|f|$, then $g\in X({\mathbb R}^n)$ with $\|g\|_X\le\|f\|_X$.

\item[{\rm (ii)}] If $0\le f_1(x)\le f_2(x)\le\cdots$ and $f_k(x)\to f(x)$ $(k\to\infty)$ hold for almost every $x\in{\mathbb R}^n$, then we have $\|f_k\|_X\to\|f\|_X$ $(k\to\infty)$.

\item[{\rm (iii)}] If a measurable set $E$ satisfies $|E|<\infty$, then we have $\|\chi_E\|_X<\infty$.

\item[{\rm (iv)}] If a measurable set $E$ satisfies $|E|<\infty$, then $\int_E|f(x)|\,{\rm d}x\le C_E\|f\|_X$ holds, where $C_E$ is a positive constant independent of $f$.
\end{itemize}
In particular, when $X({\mathbb R}^n)$ is a Banach space satisfying the four conditions above, it is called a Banach function space.
\end{definition}

The theory of Banach function spaces has been extended to ball Banach function spaces, which were defined by Sawano, Ho, Yang, and Yang in \cite{SHYY17}.

\begin{definition}\label{def:bBfs}
The (quasi-)Banach function space $X({\mathbb R}^n)$ in which the conditions (iii) and (iv) are replaced, respectively, by
\begin{itemize}
\item[${\rm (iii)}'$] For any ball $B\subset{\mathbb R}^n$, $\|\chi_B\|_X<\infty$,

\item[${\rm (iv)}'$] For any ball $B\subset{\mathbb R}^n$, $\int_B|f(x)|\,{\rm d}x\le C_B\|f\|_X$ holds, where $C_B$ is a positive constant independent of $f$,
\end{itemize}
is called a (quasi-)ball Banach function space.
\end{definition}

The conditions (i) and (ii) in Definition \ref{def:Bfs} are referred to as the ideal property and the Fatou property, respectively.

To state some statements in this paper, we use the following 2 symbols.
The associate space $X'({\mathbb R}^n)$ of $X({\mathbb R}^n)$ is defined as the space of all $f\in L^0({\mathbb R}^n)$ with the finite functional
\[
\|f\|_{X'}
\equiv
\sup_{\|g\|_X\le1}
\int_{{\mathbb R}^n}|f(x)g(x)|\,{\rm d}x.
\]
The Hardy--Littlewood maximal operator $M$ is defined by
\[
Mf(x)
\equiv
\sup_{Q\in{\mathcal Q}}
\frac{\chi_Q(x)}{|Q|}\int_Q|f(y)|\,{\rm d}y,
\quad
x\in{\mathbb R}^n,
\]
for $f\in L^0({\mathbb R}^n)$.
Then the characterizations for the Campanato norms via ball Banach function spaces are given by Izuki and Sawano in \cite{IzSa17}.

\begin{theorem}\label{thm:IzSa17}
Let $\theta\in(0,1)$, let $X({\mathbb R}^n)\subset L^0({\mathbb R}^n)$ be a Banach function space, and set $\phi_\theta:{\mathbb R}^n\times(0,\infty)\to(0,\infty)$ by $\phi_\theta(Q)=\ell(Q)^{-\theta}$.
If $M$ is bounded on $X'({\mathbb R}^n)$, then,
\[
\sup_{Q\in{\mathcal Q}}
\ell(Q)^{-\theta}
\frac{\|(f-f_Q)\chi_Q\|_X}{\|\chi_Q\|_X}
\sim
\|f\|_{{\mathcal L}_{1,\phi_\theta}}.
\]
\end{theorem}

In this paper, we aim to provide similar characterizations using Banach function spaces that do not assume the Fatou property.
It is known that, under the Fatou property, the following holds:
\[
X''({\mathbb R}^n)
=
X({\mathbb R}^n),
\]
which is known as the Lorentz--Luxemburg theorem (see Proposition \ref{prop:LN} (2)).
This allows us to use the duality argument for Banach function spaces. Moreover, the proofs of the aforementioned characterizations for the Campanato and BMO norms are also based on duality arguments.
However, there are cases where the analysis becomes difficult when the problem is transferred to the associate space via the duality argument.
For example, this occurs when the associate space cannot be explicitly identified, or even if it can be identified, when the dual space is complicated, as in the case of Morrey spaces (see \cite[Chapter 9]{SDH20}).
The aim of the present study is, therefore, to propose a method within the framework of Banach function spaces to circumvent such difficulties.

To remove the assumption of the Fatou property, we succeeded in establishing similar characterizations without using the duality argument.
More precisely, it suffices to consider Banach function spaces satisfying the only conditions (i) and ${\rm (iii)}'$ in Definitions \ref{def:Bfs} and \ref{def:bBfs}.

\begin{definition}\label{def:sBfs}
Let $X({\mathbb R}^n)=(X({\mathbb R}^n),\|\cdot\|_X)\subset L^0({\mathbb R}^n)$ be a quasi-Banach space.
The space $X({\mathbb R}^n)$ is said to be a quasi-sBfs if it satisfies the following conditions for all $f,g\in L^0({\mathbb R}^n)${\rm :}
\begin{itemize}
\item[{\rm (i)}] If $f\in X({\mathbb R}^n)$ and $g\in L^0({\mathbb R}^n)$ with $|g|\le|f|$, then $g\in X({\mathbb R}^n)$ with $\|g\|_X\le\|f\|_X$.

\item[${\rm (iii)}''$] For any cube $Q\in{\mathcal Q}({\mathbb R}^n)$, $\|\chi_Q\|_X<\infty$.
\end{itemize}
\end{definition}

Our main result is stated as follows.

\begin{theorem}\label{thm:main}
Let $X({\mathbb R}^n)$ be a quasi-sBfs, let $\phi:{\mathbb R}^n\times(0,\infty)\to(0,\infty)$ satisfy the conditions \eqref{eq:phi-dec} and \eqref{eq:near}, and let $f\in L_{\rm loc}^1({\mathbb R}^n)$.
Then the following assertions hold{\rm :}
\begin{itemize}
\item[{\rm (1)}] If there exists $C>0$ such that for all $Q\in{\mathcal Q}({\mathbb R}^n)$,
\begin{equation}\label{eq:A-X}
\|\chi_Q\|_X\|\chi_Q\|_{X'}
\le C
|Q|,
\end{equation}
then there exists $C>0$ independent of $f$ such that
\[
\|f\|_{{\mathcal L}_{1,\phi}}
\le C
\sup_{Q\in{\mathcal Q}}
\phi(Q)
\frac{\|(f-f_Q)\chi_Q\|_X}{\|\chi_Q\|_X}.
\]

\item[{\rm (2)}] If there exist $\eta>1$ and $C>0$ such that for all $\{f_j\}_{j\in{\mathbb Z}}\subset L^0({\mathbb R}^n)$,
\begin{equation}\label{eq:B-X}
\left\|
\sum_{j\in{\mathbb Z}}
(Mf_j)^\eta
\right\|_X
\le C
\left\|
\sum_{j\in{\mathbb Z}}
|f_j|^\eta
\right\|_X,
\end{equation}
then there exists $C>0$ independent of $f$ such that
\[
\sup_{Q\in{\mathcal Q}}
\phi(Q)
\frac{\|(f-f_Q)\chi_Q\|_X}{\|\chi_Q\|_X}
\le C
\|f\|_{{\mathcal L}_{1,\phi}}.
\]
\end{itemize}
\end{theorem}

Note that, when the condition ${\rm (iii)}''$ in Definition \ref{def:sBfs} is assumed, the associate norm $\|\cdot\|_{X'}$ is well defined (see Remark \ref{rem:associate}).
Moreover, under the ideal property, it is straightforward to see that conditions ${\rm (iii)}'$ and ${\rm (iii)}''$ in Definitions \ref{def:bBfs} and \ref{def:sBfs}, respectively, are equivalent.

Meanwhile, according to \cite[Lemma 2.2]{Izuki15}, if $M$ is weakly bounded on $X({\mathbb R}^n)$, that is, for all $\lambda>0$ and $f\in X({\mathbb R}^n)$, there exists a constant $C>0$ such that
\[
\lambda
\|
\chi_{\{x\in{\mathbb R}^n\,:\,Mf(x)>\lambda\}}
\|_X
\le C
\|f\|_X,
\]
then the assumption \eqref{eq:A-X} in Theorem \ref{thm:main} holds.
Of course, if we assume only the ideal property for the quasi-Banach space $X({\mathbb R}^n)$, the same statement can be established by a similar argument.
Hence, the following result is obtained.

\begin{lemma}\label{lem:Izuki15}
Let $X({\mathbb R}^n)$ be a quasi-Banach space satisfying the ideal property.
If $M$ is weakly bounded on $X({\mathbb R}^n)$, then the condition \eqref{eq:A-X} holds.
\end{lemma}

By this lemma, we can state the main theorem without using the associate functional $\|\cdot\|_{X'}$ as follows.

\begin{theorem}
Let $X({\mathbb R}^n)$ be a quasi-sBfs, let $\phi:{\mathbb R}^n\times(0,\infty)\to(0,\infty)$ satisfy the conditions \eqref{eq:phi-dec} and \eqref{eq:near}, and let $f\in L_{\rm loc}^1({\mathbb R}^n)$.
Then the following assertions hold{\rm :}
\begin{itemize}
\item[{\rm (1)}] If the Hardy--Littlewood maximal operator $M$ is weakly bounded on $X({\mathbb R}^n)$, then there exists $C>0$ independent of $f$ such that
\[
\|f\|_{{\mathcal L}_{1,\phi}}
\le C
\sup_{Q\in{\mathcal Q}}
\phi(Q)
\frac{\|(f-f_Q)\chi_Q\|_X}{\|\chi_Q\|_X}.
\]

\item[{\rm (2)}] If there exist $\eta>1$ and $C>0$ such that for all $\{f_j\}_{j\in{\mathbb Z}}\subset L^0({\mathbb R}^n)$,
\begin{equation*}
\left\|
\sum_{j\in{\mathbb Z}}
(Mf_j)^\eta
\right\|_X
\le C
\left\|
\sum_{j\in{\mathbb Z}}
|f_j|^\eta
\right\|_X,
\end{equation*}
then there exists $C>0$ independent of $f$ such that
\[
\sup_{Q\in{\mathcal Q}}
\phi(Q)
\frac{\|(f-f_Q)\chi_Q\|_X}{\|\chi_Q\|_X}
\le C
\|f\|_{{\mathcal L}_{1,\phi}}.
\]
\end{itemize}
\end{theorem}

We adopt standard notation for inequalities throughout this paper.
We use $C$ to denote a positive constant, which may vary from line to line.
If $A\le CB$, we write $A\lesssim B$ or $B\gtrsim A$.
If both $A\lesssim B$ and $A\gtrsim B$ hold, we write $A\sim B$.

The remainder of this paper is organized as follows.
Section \ref{s:Pre} presents a key lemma used in the proof of the main theorem.
In Section \ref{s:proof-main}, we prove Theorem \ref{thm:main}.
Section \ref{s:remarks} provides some detailed remarks on Banach function spaces.
Section \ref{s:examples-remarks} gives several examples of quasi-sBfs $X({\mathbb R}^n)$ appearing in Theorem \ref{thm:main}, along with comparisons to known results.
Finally, in Section \ref{s:average}, we present an alternative type of characterization, using the average quasi-norm in place of the ratio form.

\section{Preliminaries}\label{s:Pre}

To prove Theorem \ref{thm:main}, we employ the Calder\'on--Zygmund decomposition using a sparse family.
A collection ${\mathcal S}\subset{\mathcal Q}({\mathbb R}^n)$ is called a sparse family if for every $Q\in{\mathcal S}$, there exists a measurable subset $E_Q\subset Q$ such that:
\begin{itemize}
\item[(1)] $|Q|\le2|E_Q|$;
\item[(2)] The family $\{E_Q\}_{Q\in{\mathcal S}}$ is pairwise disjoint.
\end{itemize}
In addition, ${\mathcal D}(Q)$ denotes the collection of all cubes obtained by finitely many bisections of $Q$.
Using the Calder\'on--Zygmund decomposition, we obtain the following lemma, which is a simplified version of Lemma 5.1 in \cite{LOR17} (see also \cite[Lemma 3.1.2]{Hytonen21}).

\begin{lemma}\label{lem:BMO-CZ}
Let $Q_0\in{\mathcal Q}({\mathbb R}^n)$.
Then there exists a sparse family ${\mathcal S}\subset{\mathcal D}(Q_0)$ such that
\[
|f(x)-f_{Q_0}|
\lesssim
\sum_{Q\in{\mathcal S}}
\frac{\chi_Q(x)}{|Q|}
\int_Q|f(y)-f_Q|\,{\rm d}y
\quad
\text{a.e. $x\in Q_0$}.
\]
\end{lemma}

\section{Proof of Theorem \ref{thm:main}}\label{s:proof-main}

In  this section, we present the proof of Theorem \ref{thm:main}.

(1) By H\"older's inequality,
\begin{align*}
\frac{\phi(Q)}{|Q|}
\int_Q|f(x)-f_Q|\,{\rm d}x
&\le
\frac{\phi(Q)}{|Q|}
\|(f-f_Q)\chi_Q\|_X\|\chi_Q\|_{X'}\\
&\lesssim
\phi(Q)
\frac{\|(f-f_Q)\chi_Q\|_X}{\|\chi_Q\|_X}
\end{align*}
for each $Q\in{\mathcal Q}({\mathbb R}^n)$.
Then,
\[
\|f\|_{{\mathcal L}_{1,\phi}}
\lesssim
\sup_{Q\in{\mathcal Q}}
\phi(Q)
\frac{\|(f-f_Q)\chi_Q\|_X}{\|\chi_Q\|_X}.
\]

(2) Our argument in this case follows the approach of \cite{OPRR20}.
By Lemma \ref{lem:BMO-CZ}, there exists a sparse family ${\mathcal S}\subset{\mathcal D}(Q_0)$ such that
\[
\|(f-f_{Q_0})\chi_{Q_0}\|_X
\lesssim
\left\|
\sum_{Q\in{\mathcal S}}
\frac{\chi_Q}{|Q|}
\int_Q|f(y)-f_Q|\,{\rm d}y
\right\|_X.
\]
Since there exists a pairwise disjoint family $\{E_Q\}_{Q\in{\mathcal S}}$ such that
\[
E_Q\subset Q,
\quad
|Q|\le2|E_Q|,
\]
we have $\chi_Q\lesssim M\chi_{E_Q}$, and then
\begin{align*}
\|(f-f_{Q_0})\chi_{Q_0}\|_X
&\lesssim
\left\|
\sum_{Q\in{\mathcal S}}
\frac{(M\chi_{E_Q})^\eta}{|Q|}
\int_Q|f(y)-f_Q|\,{\rm d}y
\right\|_X\\
&\lesssim
\left\|
\sum_{Q\in{\mathcal S}}
\frac{\chi_{E_Q}}{|Q|}
\int_Q|f(y)-f_Q|\,{\rm d}y
\right\|_X\\
&\le
\|\chi_{Q_0}\|_X
\|f\|_{{\mathcal L}_{1,\phi}}
\sup_{Q\in{\mathcal D}(Q_0)}
\frac1{\phi(Q)},
\end{align*}
where in the second inequality we have used the assumption.
Here, by \eqref{eq:phi-dec} and \eqref{eq:near}, when $Q\subset Q_0$,
\[
\phi(Q)
\gtrsim
\phi(Q_0).
\]
Then we conclude that
\[
\phi(Q_0)
\frac{\|(f-f_{Q_0})\chi_{Q_0}\|_X}{\|\chi_{Q_0}\|}
\lesssim
\|f\|_{{\mathcal L}_{1,\phi}}.
\]
We finish the proof Theorem \ref{thm:main}.

\section{Remarks}\label{s:remarks}

To compile this section, we referred to the survey paper \cite{LoNi24} by Lorist and Nieraeth.
The original Banach function spaces not assuming the Fatou property were first considered by Zaanen and Luxemburg \cite{LuZa63,Zaanen67}.

In this section, the symbol $L^0(\Omega)$ denotes the set of all measurable functions on a general measure space $(\Omega,\mu)$.

\begin{definition}[{\cite[Section 2]{LoNi24}}]\label{def:Bfs-original}
The space $X(\Omega)=(X(\Omega),\|\cdot\|_X)\subset L^0(\Omega)$ called a quasi-Banach function space if it satisfies the
following properties.
\begin{itemize}
\item[{\rm (i)}] If $f\in X(\Omega)$ and $g\in L^0(\Omega)$ with $|g|\le|f|$, then $g\in X(\Omega)$ with $\|g\|_X\le\|f\|_X$.

\item[{\rm (v)}] For every $\mu$-measurable set $E\subset\Omega$ of positive measure, there exists a $\mu$-measurable set $F\subset E$ of positive measure with $\chi_F\in X(\Omega)$.
\end{itemize}
\end{definition}

The condition (v) in this definition is called the saturation property, and in this paper, we use a quasi-Banach function space defined by replacing this condition with the analogous condition ${\rm (iii)}''$.
More precisely, although the condition ${\rm (iii)}''$ in Definition \ref{def:sBfs} is stronger than the condition (v).
Note that the saturation property can also be defined for quasi-Banach spaces on a general measure space without assuming any metric structure.

The following equivalent formulations are provided.

\begin{proposition}\label{prop:LN}
Let $(\Omega,\mu)$ be a $\sigma$-finite measure space, and let $X(\Omega)\subset L^0(\Omega)$ be a quasi-Banach space satisfying the ideal property.
Then the following assertions hold.
\begin{itemize}
\item[{\rm (1)}] {\rm \cite[Proposition 2.5]{LoNi24}}
$X(\Omega)$ satisfies the saturation property if and only if
\[
\int_\Omega|f(x)g(x)|\,{\rm d}\mu(x)
=0
\]
for all $f\in X(\Omega)$ implies $g=0$ $\mu$-a.e.

\item[{\rm (2)}] {\rm \cite[Theorem 3.6]{LoNi24}}
Assume that $X(\Omega)$ is a Banach function space defined in Definition {\rm \ref{def:Bfs-original}}.
Then $X(\Omega)$ satisfies the Fatou property if and only if $X''(\Omega)=X(\Omega)$ with equal norm.
\end{itemize}
\end{proposition}

\begin{remark}\label{rem:associate}
It is straightforward to verify that the associate functional $\|\cdot\|_{X'}$ satisfies the triangle inequality.
Moreover, by Proposition \ref{prop:LN} (1), if $X(\Omega)$ satisfies the saturation property, then the associate functional $\|\cdot\|_{X'}$ has the property that
\[
\|f\|_{X'}=0
\quad \Longrightarrow \quad
f=0
\; \text{$\mu$-a.e.}
\]
Therefore, in particular, for the quasi-sBfs $X(\mathbb{R}^n)$, the functional $\|\cdot\|_{X'}$ is a norm.
\end{remark}

\section{Examples and comparison with the known results}\label{s:examples-remarks}

When $X''(\mathbb{R}^n)=X(\mathbb{R}^n)$, the condition \eqref{eq:A-X} is equivalent to the so-called $A_X$-condition, which was introduced in \cite{Lerner23} in connection with a generalization of Muckenhoupt's weighted theory.
We say that the $A_X$-condition holds if, for every locally integrable function $f$ and every cube $Q \in \mathcal{Q}(\mathbb{R}^n)$,
\[
|f|_Q\|\chi_Q\|_X
\lesssim
\|f\chi_Q\|_X.
\]
Remark that the $A_X$-condition follows from the weakly boundedness of $M$ on $X(\mathbb{R}^n)$ (see Lemma \ref{lem:Izuki15} above).
Hence, if $X(\mathbb{R}^n)$ is a Banach function space satisfying the Fatou property, then by a duality argument, the condition \eqref{eq:A-X} follows from the boundedness of $M$ on the associate space $X'(\mathbb{R}^n)$.

Meanwhile, the assumption \eqref{eq:B-X} in Theorem \ref{thm:main} can be regarded as an improvement of the assumption that ``$M$ is bounded on $X'(\mathbb{R}^n)$'' in Theorem \ref{thm:IzSa17}, by extending it to quasi-sBfs.
In fact, using the duality argument, we obtain
\begin{align*}
\left\|
\sum_{j\in{\mathbb Z}}
(Mf_j)^\eta
\right\|_X
&=
\sup_{\|g\|_{X'}\le1}
\sum_{j\in{\mathbb Z}}
\int_{{\mathbb R}^n}
Mf_j(x)^\eta|g(x)|
\,{\rm d}x\\
&\lesssim
\sup_{\|g\|_{X'}\le1}
\sum_{j\in{\mathbb Z}}
\int_{{\mathbb R}^n}
|f_j(x)|^\eta Mg(x)
\,{\rm d}x\\
&\le
\sup_{\|g\|_{X'}\le1}
\left\|
\sum_{j\in{\mathbb Z}}
|f_j|^\eta
\right\|_X
\|Mg\|_{X'}
\lesssim
\left\|
\sum_{j\in{\mathbb Z}}
|f_j|^\eta
\right\|_X,
\end{align*}
where the first inequality follows from the Fefferman--Stein dual inequality.

Then, in particular, the Lebesgue spaces $X({\mathbb R}^n)=L^p({\mathbb R}^n)$ $(0<p<\infty)$ satisfy the assumption \eqref{eq:B-X}.
However, the space $X({\mathbb R}^n)=L^\infty({\mathbb R}^n)$ does not satisfy the assumption \eqref{eq:B-X}.
More precisely, the estimate
\[
\sup_{Q\in\mathcal Q}
\phi(Q)
\frac{\|(f-f_Q)\chi_Q\|_{L^\infty}}{\|\chi_Q\|_{L^\infty}}
\lesssim
\|f\|_{\mathcal L_{1,\phi}}
\]
does not hold in general.
In fact, we obtain the following proposition.

\begin{proposition}
Let $\phi=\phi(\cdot):(0,\infty)\to(0,\infty)$ be a function
independent of the spatial variable $x\in{\mathbb R}^n$
satisfying \eqref{eq:phi-dec}.
Then the following assertions hold{\rm :}
\begin{itemize}
\item[{\rm (1)}] If
\[
r\phi(r)\lesssim t\phi(t),
\quad
0<r<t<\infty,
\]
then the function
\[
f(x)
=
\int_{\min(1,|x|)}^2
\frac1{\phi(t)}
\,\frac{{\rm d}t}t
\]
belongs to $\mathcal L_{1,\phi}({\mathbb R}^n)$.

\item[{\rm (2)}] If
\[
\int_0^r\frac1{\phi(t)}\,\frac{{\rm d}t}t
=
-\infty
\]
for any sufficiently small $r>0$, then for the function $f$
given in {\rm (1)},
\[
\|(f-f_{B(r)})\chi_{B(r)}\|_{L^\infty}
=
\infty,
\quad
r\in(0,1),
\]
where $B(r)$ denotes the ball centered at the origin with radius
$r>0$, and $f_{B(r)}$ is the average of $f$ over $B(r)$.
\end{itemize}
\end{proposition}

\begin{proof}
\begin{itemize}
\item[(1)]
This example is given in \cite{NaYa85}
(see also \cite{Spanne65} for the proof).

\item[(2)]
Since $r\in(0,1)$, we have
\[
f_{B(r)}
=
\frac1{r^n}
\int_0^2
\frac{\min(t^n,r^n)}{\phi(t)}
\,\frac{{\rm d}t}t.
\]
Then, for all $\rho\in(0,r)$ and any unit vector
$\omega\in{\mathbb R}^n$,
\[
f(\rho\omega)-f_{B(r)}
=
\int_\rho^r
\frac1{\phi(t)}
\,\frac{{\rm d}t}t
-
\frac1{r^n}
\int_0^r
\frac{t^n}{\phi(t)}
\,\frac{{\rm d}t}t.
\]
Here, by the almost decreasing property of $\phi$,
\[
\frac1{r^n}
\int_0^r
\frac{t^n}{\phi(t)}
\,\frac{{\rm d}t}t
\lesssim
\frac1{n\phi(r)}
<
\infty.
\]
Therefore, by the assumption,
\[
\lim_{\rho\downarrow0}
|f(\rho\omega)-f_{B(r)}|
=
\infty,
\]
which completes the proof.
\end{itemize}
\end{proof}

Here and below, we provide some examples of quasi-sBfs $X({\mathbb R}^n)$ that satisfy the conditions \eqref{eq:A-X} and \eqref{eq:B-X} in Theorem \ref{thm:main}, and compare them with those in Theorem \ref{thm:IzSa17}.

\subsection{Weighted Lebesgue spaces}

Nonnegative locally integrable functions are called weight, and for a weight $w$,
\[
w(E)
\equiv
\int_Ew(x)\,{\rm d}x
\]
is a weighted measure of a measurable set $E\subset{\mathbb R}^n$.
Then the weighted Lebesgue space $L^p({\mathbb R}^n,w)$ ($1\le p<\infty$) is defined by the space of all $f\in L^0({\mathbb R}^n)$ with the finite quasi-norm
\[
\|f\|_{L^p(w)}
\equiv
\left(
\int_{{\mathbb R}^n}
|f(x)|^pw(x)
\,{\rm d}x
\right)^{\frac1p}.
\]
Remark that $A_{L^p(w)}$-condition holds if and only if $w$ is a  Muckenhoupt's $A_p$-weight denoted by $w\in A_p$ (see \cite{Muckenhoupt72}).
Addtionally, when $w\in A_p$ for $1<p<\infty$, it is known that the vector-valued maximal inequality for the weighted Lebesgue space $L^p({\mathbb R}^n,w)$ holds.

\begin{theorem}[{\cite{AnJo80}}]
Let $1<p,\eta<\infty$, and let $w$ be a weight.
If $w\in A_p$, then for all $\{f_j\}_{j\in{\mathbb Z}}\subset L^p({\mathbb R}^n,w)$,
\[
\left\|
\left(
\sum_{j\in{\mathbb Z}}
(Mf_j)^\eta
\right)^{\frac1\eta}
\right\|_{L^p(w)}
\lesssim
\left\|
\left(
\sum_{j\in{\mathbb Z}}
|f_j|^\eta
\right)^{\frac1\eta}
\right\|_{L^p(w)}.
\]
\end{theorem}

Then, $w\in A_{p\eta}$ for $\eta>\max(1,1/p)$ implies
\[
\left\|
\sum_{j\in{\mathbb Z}}
(Mf_j)^\eta
\right\|_{L^p(w)}
\lesssim
\left\|
\sum_{j\in{\mathbb Z}}
|f_j|^\eta
\right\|_{L^p(w)}.
\]
It follows that if we assume that $w$ is an $A_\infty$-weight denoted by $w\in A_\infty$ which is defined by $w\in A_p$ for some $p\ge1$, then the weighted Lebesgue space $L^p({\mathbb R}^n,w)$ satisfies the condition \eqref{eq:B-X}.
Therefore, for the weighted Lebesgue spaces, we obtain the following result.

\begin{theorem}\label{thm:Lp-w}
Let $0<p<\infty$, let $f\in L_{\rm loc}^1({\mathbb R}^n)$, and let $w$ be a weight.
Assume that the map $\phi:{\mathbb R}^n\times(0,\infty)\to(0,\infty)$ satisfies \eqref{eq:phi-dec} and \eqref{eq:near}.
Then the following assertions hold{\rm :}
\begin{itemize}
\item[{\rm (1)}] If $1\le p<\infty$ and $w\in A_p$, then
\[
\|f\|_{{\mathcal L}_{1,\phi}}
\sim
\sup_{Q\in{\mathcal Q}}
\phi(Q)
\frac{\|(f-f_Q)\chi_Q\|_{L^p(w)}}{\|\chi_Q\|_{L^p(w)}}.
\]

\item[{\rm (2)}] If $w\in A_\infty$, then
\[
\sup_{Q\in{\mathcal Q}}
\phi(Q)
\frac{\|(f-f_Q)\chi_Q\|_{L^p(w)}}{\|\chi_Q\|_{L^p(w)}}
\lesssim
\|f\|_{{\mathcal L}_{1,\phi}}.
\]
\end{itemize}
\end{theorem}

For the case $\phi(x,r)\equiv1$, namely, ${\mathcal L}_{1,\phi}({\mathbb R}^n)={\rm BMO}({\mathbb R}^n)$, this theorem coincides with the theorem given in \cite[Section 3]{Ho11}.
Meanwhile, it is provided in \cite{MuWh75} that when $w\in A_\infty$, the norm characterization
\[
\sup_{Q\in{\mathcal Q}}
\frac1{w(Q)}
\int_Q|f(x)-f_Q|w(x)\,{\rm d}x
\sim
\|f\|_{\rm BMO}
\]
holds (see also \cite{OPRR20,Tsutsui14}).
But this inequality \lq\lq $\gtrsim$'' is not included in Theorem \ref{thm:Lp-w}.

\subsection{Lorentz and weighted Lorentz spaces}\label{ss:Lorentz}

Lorentz spaces are introduced by Lorentz \cite{Lorentz50} in 1950.
When $1<p<\infty$ and $ 0<q\le\infty$, Lorentz space $L^{p,q}({\mathbb R}^n)$, which is defined as follows, satisfies the conditions \eqref{eq:A-X} and \eqref{eq:B-X}.
On the other hand, when the quasi-Banach lattice $X({\mathbb R}^n)$ is weighted Lorentz spaces, Theorem \ref{thm:main} happens some different phenomena.

\begin{definition}
Let $0<p<\infty$ and $0<q\le\infty$, and let $w$ be a weight.
The weighted Lorentz space $L^{p,q}({\mathbb R}^n,w)$ is defined by the space of all $f\in L^0({\mathbb R}^n)$ with the finite quasi-norm
\[
\|f\|_{L^{p,q}(w)}
\equiv
\begin{cases}
\displaystyle
\left(
\int_0^\infty
\left[
t^{\frac1p}f_w^\ast(t)
\right]^q
\,\frac{{\rm d}t}t
\right)^{\frac1q},
& q<\infty, \\
\displaystyle
\sup_{t>0}
t^{\frac1p}f_w^\ast(t),
& q=\infty,
\end{cases}
\]
where $f_w^\ast(t)$ is a weighted rearrangement decreasing of $f$, that is,
\[
f_w^\ast(t)
\equiv
\inf\{
\alpha>0\,:\,
w(\{x\in{\mathbb R}^n\,:\,|f(x)|>\alpha\})
\le t
\}.
\]
\end{definition}

Remark that when $w(x)=1$, weighted Lorentz spaces are usually Lorentz spaces, and we write $(L^{p,q}({\mathbb R}^n,w),\|\cdot\|_{L^{p,q}(w)})$ by $(L^{p,q}({\mathbb R}^n),\|\cdot\|_{L^{p,q}})$.
Here, the boundedness of $M$ on the weighted Lorentz spaces are given by Chung et al. in \cite{CHK82,HuKu83}.
The class $A(p,1)$ denotes the set of all weights for which the $A_{L^{p,1}}$-condition holds.

\begin{theorem}
Let $1<p<\infty$ and $1\le q\le\infty$, and let $w$ be a weight.
\begin{itemize}
\item[{\rm (1)}] $w\in A_p$
if and only if
$M$ is bounded on $L^{p,q}(w)$.

\item[{\rm (2)}] $w\in A(p,1)$
if and only if
$M$ is weakly bounded on $L^{p,1}(w)$.
\end{itemize}
\end{theorem}

Moreover, the condition $w\in A(p,1)$ has the following statement.

\begin{lemma}[{\cite{CHK82}}]\label{lem:CHK82}
Let $1\le p<\infty$.
Then the following assertions hold{\rm :}
\begin{itemize}
\item[{\rm (1)}] $w\in A(p,1)$ if and only if for all $Q\in{\mathcal Q}({\mathbb R}^n)$ and all $E\subset Q$,
\begin{equation}\label{eq:Ap1-weight}
\frac{|E|}{|Q|}
\lesssim
\left(
\frac{w(E)}{w(Q)}
\right)^{\frac1p}.
\end{equation}
Especially, $A_{(1,1)}=A_1$, and when $1<p<\infty$, $A_p\subsetneq A(p,1)$ hold.

\item[{\rm (2)}] $\bigcup_{p\ge1}A(p,1)=A_\infty$.

\item[{\rm (3)}] If $1\le p<q<\infty$, then $A_{(p,1)}\subsetneq A_q$.
\end{itemize}
\end{lemma}

Therefore, when $X({\mathbb R}^n)=L^{p,q}({\mathbb R}^n,w)$, the following result is obtained from Theorem \ref{thm:main}.

\begin{theorem}\label{thm:Lpq(w)}
Let $1\le p<\infty$ and $1\le q\le\infty$, let $f\in L_{\rm loc}^1({\mathbb R}^n)$, and let $w$ be a weight.
Assume that the map $\phi:{\mathbb R}^n\times(0,\infty)\to(0,\infty)$ satisfies \eqref{eq:phi-dec} and \eqref{eq:near}.
Then the following assertions hold{\rm :}
\begin{itemize}
\item[{\rm (1)}] If either \lq\lq $1<p<\infty$ and $w\in A_p$'' or \lq\lq $p=q=1$ and $w\in A_1$'', then
\[
\sup_{Q\in{\mathcal Q}}
\phi(Q)
\frac{\|(f-f_Q)\chi_Q\|_{L^{p,q}(w)}}{\|\chi_Q\|_{L^{p,q}(w)}}
\sim
\|f\|_{{\mathcal L}_{1,\phi}}.
\]

\item[{\rm (2)}] If $w\in A(p,1)$, then
\[
\sup_{Q\in{\mathcal Q}}
\phi(Q)
\frac{\|(f-f_Q)\chi_Q\|_{L^{p,1}(w)}}{\|\chi_Q\|_{L^{p,1}(w)}}
\sim
\|f\|_{{\mathcal L}_{1,\phi}}.
\]
\end{itemize}
\end{theorem}

According to Lemma \ref{lem:CHK82} (1), in Theorem \ref{thm:Lpq(w)}, the assumption of (1) is included in the assumption of (2).
Additionally, the case $p=q=1$ in Theorem \ref{thm:Lpq(w)} (1) and the case $p=1$ in Theorem \ref{thm:Lpq(w)} (2) are equivalent to the case $p=1$ in Theorem \ref{thm:Lp-w} (1).

\subsection{Orlicz spaces}

Orlicz spaces are introduced by Birnbaum and Orlicz \cite{BiOr31} in 1931.
The mapping $\Phi:[0,\infty)\to[0,\infty)$ is called Young function when it satisfies the following conditions:
\begin{itemize}
\item[(1)] $\Phi$ is a positive function on $(0,\infty)$.
\item[(2)] $\Phi$ is a convex function.
\item[(3)] $\lim\limits_{t\downarrow0}\Phi(t)=\Phi(0)=0$.
\end{itemize}

Then the Orlicz space is defined as follows.

\begin{definition}\label{def:Orlicz-space}
Let $\Phi:[0,\infty)\to[0,\infty)$ be a Young function.
The Orlicz space $L^\Phi({\mathbb R}^n)$ is defined by
\[
L^\Phi({\mathbb R}^n)
\equiv
\left\{
f\in L^0({\mathbb R}^n)
\,:\,
\text{
$\displaystyle
\int_{{\mathbb R}^n}\Phi(k|f(x)|)\,{\rm d}x
<\infty
$,
for every $k>0$
}
\right\}
\]
endowed with the norm
\[
\|f\|_{L^\Phi}
\equiv
\inf\left\{
\lambda>0\,:\,
\int_{{\mathbb R}^n}
\Phi\left(\frac{|f(x)|}\lambda\right)
\,{\rm d}x
\le1
\right\}.
\]
\end{definition}

When $\Phi(t)=t^p$ ($1\le p<\infty$), $\|f\|_{L^\Phi}=\|f\|_{L^p}$ holds.
Thus the Orlicz spaces are generalization for the Lebesgue spaces.
Here, we consider the following conditions for the Young functions:

\begin{itemize}
\item A Young function $\Phi:[0,\infty)\to[0,\infty)$ is said to satisfy the $\Delta_2$-condition or the doubling condition, denoted by $\Phi\in\Delta_2$, 
if there exists a constant $k >1$ such that
\begin{equation*}
\Phi(2r)\le k\Phi(r)
\quad \text{for} \quad
r>0.
\end{equation*}

\item A Young function $\Phi:[0,\infty)\to[0,\infty)$ is said to satisfy the $\nabla_2$-condition, denoted by $\Phi\in\nabla_2$, 
if there exists a constant $k >1$ such that
\begin{equation*}
\Phi(r)\le\frac1{2k}\Phi(kr)
\quad \text{for} \quad
r>0.
\end{equation*}
\end{itemize}

Then the vector-valued maximal inequality for the Orlicz spaces are given as follows.

\begin{theorem}
Let $\Phi:[0,\infty)\to[0,\infty)$ be a Young function, and let $1<\eta<\infty$.
Then the following assertions hold{\rm :}
\begin{itemize}
\item[{\rm (1)}] {\rm \cite[Lemma 1.2.4]{KoKr91}}
$M$ is weakly bounded on $L^\Phi({\mathbb R}^n)$.

\item[{\rm (2)}] {\rm \cite[Theorem 2.6]{NaSa14}}
If $\Phi\in\Delta_2\cap\nabla_2$, then for all $\{f_j\}_{j\in{\mathbb Z}}\subset L^\Phi({\mathbb R}^n)$,
\[
\left\|
\left(
\sum_{j\in{\mathbb Z}}
(Mf_j)^\eta
\right)^{\frac1\eta}
\right\|_{L^\Phi}
\lesssim
\left\|
\left(
\sum_{j\in{\mathbb Z}}
|f_j|^\eta
\right)^{\frac1\eta}
\right\|_{L^\Phi}.
\]
\end{itemize}
\end{theorem}

By Lemma \ref{lem:Izuki15} and (1) in this theorem, we can check the condition of Theorem \ref{thm:main} (1).
Moreover, by (2) in this theorem, we can check the condition of Theorem \ref{thm:main} (2).
For $\theta\ge1$, we set
\[
\Phi_{\theta}(t)
\equiv
\int_0^{t^{\theta}}\frac{\Phi(s)}s\,{\rm d}s,
\quad
t\ge0.
\]
Then $\theta>1$ implies $\Phi_{\theta}\in\Delta_2\cap\nabla_2$
according to \cite{HKO23}, and we obtain
\begin{align*}
\left\|
\sum_{j=1}^\infty(M f_j)^\eta
\right\|_{L^{\Phi_1}}^{\frac1\eta}
&=
\left\|
\left(\sum_{j=1}^\infty(M f_j)^\eta\right)^{\frac1\eta}
\right\|_{L^{\Phi_\eta}}
\lesssim
\left\|
\left(\sum_{j=1}^\infty|f_j|^\eta\right)^{\frac1\eta}
\right\|_{L^{\Phi_\eta}}\\
&=
\left\|
\sum_{j=1}^\infty |f_j|^\eta
\right\|_{L^{\Phi_1}}^{\frac1\eta}.
\end{align*}
Moreover, $\|\cdot\|_{L^{\Phi_1}}$ and $\|\cdot\|_{L^{\Phi}}$ are equivalent according to \cite{HKO23}.
Thus, assuming the condition $\Phi\in\Delta_2$, we can see that the Orlicz space $X({\mathbb R}^n)=L^\Phi({\mathbb R}^n)$ satisfies the assumption of Theorem \ref{thm:main} (2).

Therefore, for the Orlicz spaces, we have the following result by Theorem \ref{thm:main}.

\begin{theorem}
Let $\Phi:[0,\infty)\to[0,\infty)$ be a Young function, let $f\in L_{\rm loc}^1({\mathbb R}^n)$, and let $\phi:{\mathbb R}^n\times(0,\infty)\to(0,\infty)$ satisfy \eqref{eq:phi-dec} and \eqref{eq:near}.
Then the following assertions hold{\rm :}
\begin{itemize}
\item[{\rm (1)}] $\displaystyle
\|f\|_{{\mathcal L}_{1,\phi}}
\lesssim
\sup_{Q\in{\mathcal Q}}
\phi(Q)
\frac{\|(f-f_Q)\chi_Q\|_{L^\Phi}}{\|\chi_Q\|_{L^\Phi}}
$.

\item[{\rm (2)}] If $\Phi\in\Delta_2$, then
\[
\sup_{Q\in{\mathcal Q}}
\phi(Q)
\frac{\|(f-f_Q)\chi_Q\|_{L^\Phi}}{\|\chi_Q\|_{L^\Phi}}
\lesssim
\|f\|_{{\mathcal L}_{1,\phi}}.
\]
\end{itemize}
\end{theorem}

\subsection{Variable Lebesgue spaces}

For a measurable map $p(\cdot):{\mathbb R}^n\to(0,\infty)$, we define
\[
p_-
\equiv
\inf_{x\in{\mathbb R}^n}
p(x),
\quad
p_+
\equiv
\sup_{x\in{\mathbb R}^n}
p(x).
\]

\begin{definition}
Let $p(\cdot):{\mathbb R}^n\to(0,\infty)$ be a measurable map.
The variable Lebesgue space $L^{p(\cdot)}({\mathbb R}^n)$ is defined by
\[
L^{p(\cdot)}({\mathbb R}^n)
\equiv
\{
f\in L^0({\mathbb R}^n)
\,:\,
\text{
$\displaystyle
\int_{{\mathbb R}^n}(k|f(x)|)^{p(x)}\,{\rm d}x
<\infty
$ for some $k>0$
}
\}
\]
endowed with the quasi-norm
\[
\|f\|_{L^{p(\cdot)}}
\equiv
\inf\left\{
\lambda>0
\,:\,
\int_{{\mathbb R}^n}\left(\frac{|f(x)|}\lambda\right)^{p(x)}\,{\rm d}x
\le1
\right\}.
\]
\end{definition}

It is not sufficient condition of the boundedness of $M$ on the variable Lebesgue space $L^{p(\cdot)}({\mathbb R}^n)$ that $A_{L^{p(\cdot)}}$-condition holds (see \cite[Theorem 5.3.4]{DHHR17}).
To state the statements for the boundedness of $M$ on the variable Lebesgue spaces and their vector-valued extenstions, we introduce the following log-H\"older conditions:
\begin{itemize}
\item[(1)] We say that $p(\cdot)$ is locally log-H\"older continuous, and denote this by $p(\cdot)\in LH_0$, if for all $x,y\in{\mathbb R}^n$ with $|x-y|<1/2$,
\[
|p(x)-p(y)|
\lesssim
\frac1{-\log|x-y|}.
\]

\item[(2)] We say that $p(\cdot)$ is log-H\"older continuous at infinity, and denote this by $p(\cdot)\in LH_\infty$, if there exists $p_\infty>0$ such that for all $x,y\in{\mathbb R}^n$,
\[
|p(x)-p_\infty|
\lesssim
\frac1{\log(e+|x|)}.
\]
\end{itemize}

\begin{proposition}
Let $p(\cdot):{\mathbb R}^n\to(0,\infty)$ be a measurable map satisfying $0<p_-\le p_+<\infty$ and $p(\cdot)\in LH_0\cap LH_\infty$, and let $1<\eta<\infty$.
Then the following assertions hold{\rm :}
\begin{itemize}
\item[{\rm (1)}] {\rm \cite[Theorem 3.16]{CrFi13}}
If $p_-\ge1$, then for all $\lambda>0$ and $f\in L^{p(\cdot)}({\mathbb R}^n)$,
\[
\lambda
\|
\chi_{\{x\in{\mathbb R}^n\,:\,Mf(x)>\lambda\}}
\|_{L^{p(\cdot)}}
\lesssim
\|f\|_{L^{p(\cdot)}}.
\]

\item[{\rm (2)}] {\rm \cite[Theorem 3.16]{CrFi13}}
If $p_->1$, then for all $f\in L^{p(\cdot)}({\mathbb R}^n)$,
\[
\|Mf\|_{L^{p(\cdot)}}
\lesssim
\|f\|_{L^{p(\cdot)}}.
\]
\end{itemize}
\end{proposition}

\begin{proposition}[{\cite{CFMP06}}]
Let $p(\cdot):{\mathbb R}^n\to(0,\infty)$ be a measurable map satisfying $1<p_-\le p_+<\infty$.
If $M$ is bounded on $L^{p(\cdot)}({\mathbb R}^n)$, then for all $\{f_j\}_{j\in{\mathbb Z}}\subset L^{p(\cdot)}({\mathbb R}^n)$,
\[
\left\|
\left(
\sum_{j\in{\mathbb Z}}
(Mf_j)^\eta
\right)^{\frac1\eta}
\right\|_{L^{p(\cdot)}}
\lesssim
\left\|
\left(
\sum_{j\in{\mathbb Z}}
|f_j|^\eta
\right)^{\frac1\eta}
\right\|_{L^{p(\cdot)}}.
\]
\end{proposition}

Then we have the following result from Theorem \ref{thm:main} and Lemma \ref{lem:Izuki15}.

\begin{theorem}
Let $p(\cdot):{\mathbb R}^n\to(0,\infty)$ be a measurable map satisfying $0<p_-\le p_+<\infty$, let $f\in L_{\rm loc}^1({\mathbb R}^n)$, and let $\phi:{\mathbb R}^n\times(0,\infty)\to(0,\infty)$ satisfy \eqref{eq:phi-dec} and \eqref{eq:near}.
Then the following assertions hold{\rm :}
\begin{itemize}
\item[{\rm (1)}] If $p_-\ge1$, then
\[
\|f\|_{{\mathcal L}_{1,\phi}}
\lesssim
\sup_{Q\in{\mathcal Q}}
\phi(Q)
\frac{\|(f-f_Q)\chi_Q\|_{L^{p(\cdot)}}}{\|\chi_Q\|_{L^{p(\cdot)}}}.
\]

\item[{\rm (2)}] If there exists $\eta>1$ with $(\eta p(\cdot))_->1$ such that $M$ is bounded on $L^{\eta p(\cdot)}({\mathbb R}^n)$, then
\[
\sup_{Q\in{\mathcal Q}}
\phi(Q)
\frac{\|(f-f_Q)\chi_Q\|_{L^{p(\cdot)}}}{\|\chi_Q\|_{L^{p(\cdot)}}}
\lesssim
\|f\|_{{\mathcal L}_{1,\phi}}.
\]
\end{itemize}
\end{theorem}

Although it is given in \cite{IST14} that if the sufficient condition of
\[
\sup_{Q\in{\mathcal Q}}
\frac{\|(f-f_Q)\chi_Q\|_{L^{p(\cdot)}}}{\|\chi_Q\|_{L^{p(\cdot)}}}
\sim
\|f\|_{\rm BMO}
\]
is the weak-type $L^{p(\cdot)}({\mathbb R}^n)$-boundedness of $M$ for the case $\phi(x,r)\equiv1$, we cannot compare with this statement.

\subsection{Morrey spaces and their generalizations}

Morrey spaces are introduced by Morrey \cite{Morrey38} in 1938.

\begin{definition}
Let $0<q\le p<\infty$.
The Morrey space ${\mathcal M}^p_q({\mathbb R}^n)$ is defined by the space of all $f\in L^0({\mathbb R}^n)$ with the finite quasi-norm
\[
\|f\|_{{\mathcal M}^p_q}
\equiv
\sup_{Q\in{\mathcal Q}}
|Q|^{\frac1p-\frac1q}
\left(
\int_Q|f(x)|^q\,{\rm d}x
\right)^{\frac1q}.
\]
\end{definition}

The weakly boundedness of $M$ and the vector-valued maximal inequality on the Morrey spaces are given as follows.

\begin{proposition}
Let $1\le q\le p<\infty$ and $1<\eta<\infty$.
Then the following assertions hold{\rm :}
\begin{itemize}
\item[{\rm (1)}] {\rm \cite{ChFr87}}
For all $\lambda>0$ and $f\in{\mathcal M}^p_q({\mathbb R}^n)$,
\[
\lambda
\|
\chi_{\{x\in{\mathbb R}^n\,:\,Mf(x)>\lambda\}}
\|_{{\mathcal M}^p_q}
\lesssim
\|f\|_{{\mathcal M}^p_q}.
\]

\item[{\rm (2)}] {\rm \cite{SaTa05,TaXu05}}
If $q>1$, then for all $\{f_j\}_{j\in{\mathbb Z}}\subset{\mathcal M}^p_q({\mathbb R}^n)$,
\[
\left\|
\left(
\sum_{j\in{\mathbb Z}}
(Mf_j)^\eta
\right)^{\frac1\eta}
\right\|_{{\mathcal M}^p_q}
\lesssim
\left\|
\left(
\sum_{j\in{\mathbb Z}}
|f_j|^\eta
\right)^{\frac1\eta}
\right\|_{{\mathcal M}^p_q}
\]
\end{itemize}
\end{proposition}

Then, by Lemma \ref{lem:Izuki15}, Theorem \ref{thm:main} gives the following result for the case $X({\mathbb R}^n)={\mathcal M}^p_q({\mathbb R}^n)$.

\begin{theorem}
Let $1\le q\le p<\infty$, let $\phi:{\mathbb R}^n\times(0,\infty)\to(0,\infty)$ satisfy the conditions \eqref{eq:phi-dec} and \eqref{eq:near}, and let $f\in L_{\rm loc}^1({\mathbb R}^n)$.
Then
\[
\sup_{Q\in{\mathcal Q}}
\phi(Q)
\frac{\|(f-f_Q)\chi_Q\|_{{\mathcal M}^p_q}}
{\|\chi_Q\|_{{\mathcal M}^p_q}}
\sim
\|f\|_{{\mathcal L}_{1,\phi}}.
\]
\end{theorem}

Moreover, we can consider many kinds of generalizations for the Morrey spaces.
The vector-valued maximal inequalities for these function spaces are given as follows:
Generalized Morrey spaces \cite{AGNS17}, Orlicz--Morrey spaces of three kinds \cite{Ho13,SHG15,GHSN16}, mixed Morrey spaces \cite{Nogayama19}, Morrey--Lorentz spaces \cite{Hatano20}, Bourgain--Morrey spaces \cite{HNSH23}.

\section{Another type of characterizations}\label{s:average}

As another type of characterization, we rewrite the $L^p$-average $(|f|^p)_Q^{1/p}$ as
\begin{align*}
\frac{\|f\chi_Q\|_{L^p}}{\|\chi_Q\|_{L^p}}
=
\left(
\frac1{|Q|}
\int_Q|f(x)|^p\,{\rm d}x
\right)^{\frac1p}
=
\|\delta^{\ell(Q)}[f\chi_Q]\|_{L^p},
\end{align*}
where $\delta^t$ with $t>0$ is the dilatation operator defined for all $f\in L^0({\mathbb R}^n)$ by $\delta^tf(x)=f(tx)$, and $\ell(Q)$ denotes the side length of $Q\in{\mathcal Q}({\mathbb R}^n)$.
We then consider replacing $L^p({\mathbb R}^n)$ with a more general function space $X({\mathbb R}^n)$ as follows (see \cite{Perez95}):
\[
\|f\|_{X,Q}
\equiv
\|\delta^{\ell(Q)}[f\chi_Q]\|_X.
\]
Furthermore, we establish the following result.

\begin{theorem}
Let $X({\mathbb R}^n)\subset L^0({\mathbb R}^n)$ be a quasi-sBfs satisfy
\[
\sup_{Q\in{\mathcal Q}}
\|1\|_{X,Q}
<\infty,
\]
let $\phi:{\mathbb R}^n\times(0,\infty)\to(0,\infty)$ satisfy the conditions \eqref{eq:phi-dec} and \eqref{eq:near}, and let $f\in L_{\rm loc}^1({\mathbb R}^n)$.
If there exists $\eta>1$ such that
\[
\left\|
\sum_{j\in{\mathbb Z}}
(Mf_j)^\eta
\right\|_X
\lesssim
\left\|
\sum_{j\in{\mathbb Z}}
|f_j|^\eta
\right\|_X,
\]
then
\[
\sup_{Q\in{\mathcal Q}}
\phi(Q)
\|f-f_Q\|_{X,Q}
\lesssim
\|f\|_{{\mathcal L}_{1,\phi}}.
\]
\end{theorem}

Since the fundamental identity $\delta^t\circ M=M\circ\delta^t$ holds, this theorem can be proven using a similar argument as in the proof of Theorem \ref{thm:main} (2).
We omit the claim concerning statement (1) in Theorem \ref{thm:main},
since we would need to assume the technical condition
\begin{equation}\label{eq:technical}
\frac1{|Q|}\int_Q|f(x)|\,{\rm d}x
\lesssim
\|f\|_{X,Q}.
\end{equation}

A typical example of $\|\cdot\|_{X,Q}$ is the Orlicz average.
It is easy to see that, for a Young function $\Phi$, the Orlicz average
\[
\|f\|_{\Phi,Q}
\equiv
\inf\left\{
\lambda>0
\,:\,
\frac1{|Q|}
\int_Q\Phi\left(\frac{|f(x)|}\lambda\right)\,{\rm d}x
\le1
\right\}
\]
coincides with $\|f\|_{L^\Phi,Q}$.
It is straightforward to verify, using Hölder's inequality and the fact that $\|1\|_{\Phi,Q} = \Phi^{-1}(1) < \infty$, that the technical condition \eqref{eq:technical} always holds for the Orlicz average.
Then, we obtain the following corollary.

\begin{corollary}
Let $\Phi:[0,\infty)\to[0,\infty)$ be a Young function, and let $\phi:{\mathbb R}^n\times(0,\infty)\to(0,\infty)$ satisfy the conditions \eqref{eq:phi-dec} and \eqref{eq:near}, and let $f\in L_{\rm loc}^1({\mathbb R}^n)$.
If $\Phi\in\Delta_2$, then
\[
\sup_{Q\in{\mathcal Q}}
\phi(Q)
\|f-f_Q\|_{\Phi,Q}
\sim
\|f\|_{{\mathcal L}_{1,\phi}}.
\]
\end{corollary}

{\bf Acknowledgements.} 
The author was supported by Grant-in-Aid for Research Activity Start-up Grant Number 23K19013 and the Grant-in-Aid for JSPS Fellows (No. 25KJ0222).

\end{document}